\documentclass{amsart}
\usepackage{amssymb}

\newcommand{\cA}{{\mathcal A}}
\newcommand{\cO}{{\mathcal O}}

\newcommand{\be}{{\bf 1}}
\newcommand{\bo}{{\bf {\bar 1}}}

\newcommand{\BZ}{{\mathbb Z}}

\newcommand{\BQ}{{\mathbb Q}}

\newcommand{\eps}{{\varepsilon}}
\newcommand{\sq}{$\square$}

\newcommand{\Ve}{\mbox{Vec}}
\newcommand{\Rep}{\mbox{Rep}}

\newcommand{\Hom}{\mbox{Hom}}

\newcommand{\Id}{\mbox{Id}}

\newcommand{\End}{\mbox{End}}

\title{Tensor categories (after P.~Deligne)}
\author{Viktor Ostrik}
\email{ostriu@math.ias.edu}
\address{Institute for Advanced Study, Einstein Drive, Princeton, NJ 08540}
\thanks{The author was partially supported by NSF grant DMS-0098830}
\date{May 2002}
\begin{document}
\begin{abstract} This is my talk at the MIT Lie Groups Seminar.
I give an exposition of a recent paper by P.~Deligne
``Cat\'egories tensorielles''.
\end{abstract} 
\maketitle
\section{Introduction}
In a recent preprint \cite{D} P.~Deligne proved that any tensor (=rigid 
symmetric monoidal abelian) category over an algebraically closed field of
characteristic 0 satisfying certain very mild conditions
comes from representations of an affine super group (see exact statement 
below). This result is interesting by itself and also has applications
in the theory of Hopf algebras (see \cite{EG1} and
especially \cite{EG2}) since it allows to classify completely for example
finite-dimensional triangular Hopf algebras.

In this note I give an exposition of Deligne's proof oriented on 
representation theorists. So I tried to be as elementary as possible.
It is assumed that the reader knows basic notions of the category theory
and is familiar with representations of the symmetric groups over the
complex numbers. This paper does not contain anything original 
(except, possibly, mistakes) and can not
be considered as a substitute for Deligne's paper. In many cases I leave
proofs (and even definitions) to the reader, in all such cases the reader is
referred to Deligne's exellent exposition. 
  
\section{Main Theorem}
\subsection{Tensor categories}
Let $k$ be an algebraically closed field of characteristic 0. 

{\bf Definition.} A tensor category $\cA$ is a small abelian $k-$linear
category endowed with a biexact $k-$linear functor $\otimes :\cA \times \cA \to
\cA$, associativity, commutativity constraints, unit object $\be$ such that
$(\cA, \otimes )$ is a symmetric monoidal category (see e.g. \cite{Mc}).
Furthermore, the category $\cA$ is assumed to be rigid, that is for any
object $X\in \cA$ there exists object $X^\vee$ and morphisms
$coev: \be \to X\otimes X^\vee$ and $ev: X^\vee \otimes X\to \be$ such that
the compositions $X\to X\otimes X^\vee \otimes X\to X$ and $X^\vee \to
X^\vee \otimes X\otimes X^\vee$ are identity morphisms. Finally it is assumed
that $\End(\be)=k$ (or equivalently $\be$ is a simple object of $\cA$, see
\cite{D1}).

One says that a category $\cA$ is finitely $\otimes-$generated if there is
an object $X\in \cA$ such that any object of $\cA$ is isomorphic to a
subquotient of an object which is a direct sum of objects $X^{\otimes n}$. 
Such an object $X$ is called a $\otimes-$generator of $\cA$.

{\bf Examples.} (0) The category $\Ve$ of finite dimensional vector spaces
over $k$ is a tensor category. This category is obviously finitely 
$\otimes-$generated.

(1) The category $\Rep(G)$ of finite dimensional representations of an affine 
algebraic group $G$ over $k$. Any faithful representation $X$ of $G$ is 
$\otimes-$generator of this category.

(2) The category $\Rep(G)$ of finite dimensional discrete (= factorizing 
through a finite quotient) representations of a profinite group $G$. This 
category is finitely $\otimes-$generated if and only if $G$ is finite.

(3) The category $s\Ve$ of finite dimensional super vector spaces. Recall
that the objects of $s\Ve$ are $\BZ/2\BZ-$graded spaces $V_0\oplus V_1$
and commutativity morphism is given by $x\otimes y\mapsto (-1)^{deg(x)deg(y)}
y\otimes x$ where $x\in V_{deg(x)}$ and $y\in V_{deg(y)}$.

(4) Let $\cO (G)$ be a supercommutative super Hopf algebra finitely generated
as an algebra (e.g. the exterior algebra $\wedge^\bullet (V)$ of a vector space
$V$ with comultiplication $x\mapsto x\otimes 1+1\otimes x$ is a super Hopf
algebra, but not a usual Hopf algebra). One considers $\cO (G)$ as the 
functions algebra on an affine  ``super group'' $G$. Then the category 
$\Rep(G)$ of finite dimensional (super) comodules over $\cO (G)$ is a finitely 
$\otimes-$generated tensor category.

(5) Let $G$ be a super group and let $\eps \in G(k)$ (that is $\eps$ is a
homomorphism $\cO (G)\to k$) such that $\eps^2=1$ (that is the map
$\cO (G)\to k$ $h\mapsto \eps(h_{(1)})\eps(h_{(2)})$ coincides with the 
counit $\epsilon$; here
$h\mapsto h_{(1)}\otimes h_{(2)}$ is the comultiplication) and inner 
automorphism of $G$ induced by $\eps$ is just the parity automorphism (that is
automorphism of $\cO (G)$ $h\mapsto \eps(h_{(1)})h_{(2)}\eps(S(h_{(3)}))$ 
coincides with $h\mapsto (-1)^{deg(h)}h$; here $h\mapsto h_{(1)}\otimes 
h_{(2)}\otimes h_{(3)}$ is twice iterated
comultiplication and $h\mapsto S(h)$ is the antipode). Consider the category
$\Rep(G,\eps)$ consisting of objects $V$ of $\Rep(G)$ such that $\eps$ acts on
$V$ by the parity automorphism. The category $\Rep(G)$ is a tensor category.
It is finitely $\otimes-$generated if and only if $\cO (G)$ is finitely 
generated algebra.

{\bf Exercise.} Show that examples (0), (1), (3), (4) above are special cases 
of example (5). (Hint: for example (4) consider semidirect product of $G$ and 
of $\BZ/2\BZ$ where $\BZ/2\BZ$ acts on $G$ via parity automorphism).

\subsection{Schur functors} For any object $X\in \cA$ we have a natural action
of the symmetric group $S_n$ on the object $X^{\otimes n}$ induced by 
commutativity isomorphisms. For any partition $\lambda$ of $n$ let $V_\lambda$
denote the corresponding irreducible representation of $S_n$. We have the
diagonal action of $S_n$ on $V_\lambda \otimes X^{\otimes n}$. 

{\bf Definition.} The {\em Schur functor} $S_\lambda$ is defined by the formula
$$S_\lambda (X)=(V_\lambda \otimes X^{\otimes n})^{S_n}:=(\sum_{g\in
S_n}g))(V_\lambda \otimes X^{\otimes n}).$$

For example if $\lambda =n$ then $S_\lambda (X)=Sym^n(X)$ is just symmetric
power and if $\lambda =1+1+\cdots +1$ then $S_\lambda (X)=\wedge^n(X)$ is
the exterior power. 

{\bf Exercise.} Let $V=V_0\oplus V_1$ be a superspace of dimension $p|q$, that
is $\dim V_0=p$ and $\dim V_1=q$. Prove that $S_\lambda (V)=0$ if and only if
there is $i>p$ such that $\lambda_i>q$ (in other words $S_\lambda (V)\ne 0$
if and only if the Young diagram of $\lambda$ lies in the union of two strips
of width $p$ in ``$\lambda_i-$direction'' and of width $q$ in 
``$i-$direction''). Hint: One has the following formula:
$$S_\lambda (X\oplus Y)=\bigoplus_{\mu,\nu}(S_\mu (X)\otimes S_\nu(Y))^{a^
\lambda_{\mu,\nu}}$$
where the summation runs over partitions $\mu, \nu$ such that $|\mu|+|\nu|=
|\lambda|$ and $a^\lambda_{\mu,\nu}$ are Littlewood-Richardson coefficients
(additional exercise: state a similar formula for $S_\lambda(X\otimes Y)$).

In particular for any object $X$ of the category $\Rep(G,\eps)$ there exists
$\lambda$ such that $S_\lambda (X)=0$. 

\subsection{Main Theorem} Here is the main result:

{\bf Theorem.} Let $\cA$ be a finitely $\otimes-$generated tensor category
such that for any $X\in \cA$ there is $\lambda$ with $S_\lambda (X)=0$.
Then $\cA$ is equivalent as tensor category to $\Rep(G,\eps)$ for some
supergroup $G$. 

{\bf Remarks.} (i) The set of objects $X\in \cA$ annihilated by some
(depending on $X$) Schur functor is stable under direct sums, tensor products,
taking dual, subquotients, extensions; any such object has finite length.
We leave this as an exercise to the reader.

 (ii) The condition that an object $X$ is annihilated by some
Schur functor is equivalent to the existence of $N$ such that 
$length(X^{\otimes n})\le N^n$ for all $n\ge 0$ (in one direction this is
a consequence of the Theorem and in other direction the decomposition
$X^{\otimes n}=\oplus_\lambda V_\lambda \otimes S_\lambda (X)$ shows that
if $S_\lambda (X)\ne 0$ for all $\lambda$ then $length(X^{\otimes n})\ge 
\sum_\lambda \dim 
V_\lambda \ge(\sum_\lambda (\dim V_\lambda )^2)^{1/2}=
\sqrt{n!}$). This condition is automatically satisfied if category
$\cA$ has only finitely many simple objects (for example if $\cA$ is the
category of representations of a finite dimensional (weak) Hopf algebra). 

As an immediate consequence one gets

{\bf Corollary.}
Assume that $\cA$ is semisimple with finitely many simple
objects. Then $\cA$ is equivalent to category $\Rep(G,\eps)$ where $G$ is
a finite group and $\eps \in G$ is a central element of order at most 2.

\subsection{Strategy of the proof} We begin with the following 

{\bf Definition.} Let $\cA$ and $\cA'$ be two tensor categories. A {\em tensor
functor} from $\cA$ to $\cA'$ is a functor $F:\cA \to \cA'$ endowed with an
isomorphism $\be \to F(\be)$ and functorial isomorphisms $F(X)\otimes F(Y)\to
F(X\otimes Y)$ compatible with associativity, commutativity and unit 
constraints.

The main difficulty of Deligne's Theorem
is the following: the tensor category $\Rep(G,\eps)$ has an additional 
structure,
the super fiber functor (that is tensor functor to the category $s\Ve$). 
Conversely it is not very hard (and is a standard exercise in Tannakian 
formalism) to prove that if a category $\cA$ admits a super fiber functor
then it is equivalent to the category $\Rep(G,\eps)$, see \cite{D1} for the
case of usual (not super) fiber functors. So we are reduced to
showing that the category $\cA$ admits a super fiber functor. For this one
generalizes the notion of super fiber functor in the following way: for a
supercommutative algebra $R$ one defines $R-$fiber functor to be exact tensor
functor to the category of $R-$modules with tensor product over $R$ as a
tensor product (note that strictly speaking category of $R-$modules is not a 
tensor category in our sense
since $\End(\be)=R\ne k$). Then one shows (this is a key result) that category
$\cA$ admits $R-$fiber functor for some (very big) algebra $R$. Then using
standard technique from algebraic geometry one deduces that category $\cA$
admits super fiber functor over $k$. 

\subsection{Some counterexamples} It is not trivial to construct an example
of tensor category with object $V$ such that $S_\lambda (V)\ne 0$ for any
$\lambda$. Here we present two such examples.

{\bf Orthosymplectic example.} (\cite{De2}) Let $t$ be an indeterminate. 
Consider the following category $OSP_{\BQ (t)}$:

Objects: finite sets; tensor product: disjoint union;
morphisms: $\Hom (X,Y)$ is free $\BQ (t)-$module
generated by bordisms from $X$ to $Y$ (= 1-dimensional manifolds with
boundary $X\sqcup Y$) modulo the relation
$ [\mbox{bordism} \sqcup \mbox{circle}]=t[\mbox{bordism}]$; the 
composition of morphisms is induced by the composition of bordisms.

Now we define the category $OSP_t$ to be the Karoubian envelope of the
category $OSP_{\BQ (t)}$. 
 
{\bf The category $GL_t$.} (\cite{De2, D1}) This example is completely
analogous to the previous example except that we consider oriented finite
sets (that is finite sets $X$ together with map $\eps :X \to \{ \pm 1\}$)
and oriented bordisms. The resulting category is denoted $GL_t$.

The categories $GL_t$ and $OSP_t$ are abelian semisimple categories; the
simple objects of these categories are absolutely simple.

{\bf Exercise.} Let $V$ be the object of $GL_t$ or $OSP_t$ corresponding
to a finite set with one element. Show that $S_\lambda (V)\ne 0$ for any
$\lambda$.

\section{Existence of a super fiber functor}
The main point of Deligne's proof is the possibility to imitate affine 
algebraic
geometry in the category $\cA$. In other words the notion of an algebra (in 
what follows ``algebra'' means a nonzero associative
commutative algebra with unit) makes sense in the category $\cA$: an algebra in
$\cA$ is an ind-object $A$ of $\cA$ endowed with multiplication morphism
$A\otimes A\to A$ satisfying certain axioms (we leave as an exercise for the
reader to state precisely these axioms). For example for any object $X\in \cA$
one defines the symmetric algebra $Sym^*(X)$ (again we leave details to the
reader). For an algebra $A$ one easily defines the notions of $A-$modules, 
homomorphisms and tensor products over $A$.

\subsection{Key Lemma} Let $A$ be an algebra in $\cA$. The notion of a rigid
$A-$module is defined exactly as before. Note that a direct summand of a
rigid $A-$module is rigid. Also symmetric powers of modules
over $A$ are defined. 

An $A-$algebra $B$ is an algebra $B$ in $\cA$ together with a homomorphism
$A\to B$. For any $A-$module $M$ one defines its {\em extension of scalars} 
to be $B-$module $M_B:=M\otimes_AB$. In particular $\be_B$ is $B$ itself 
considered as $B-$module. Clearly, extension of scalars is a functor. It is
obvious that extension of scalars of a rigid module is again rigid module.
The Schur functors are defined over an algebra $A$ and commute with extension
of scalars. 

{\bf Key Lemma.} Let $M$ be a rigid $A-$module. The existence of $A-$algebra 
$B$ such that $M_B$ has $\be_B$ as a direct summand is equivalent to the
condition $Sym^n_A(M)\ne 0$ for all $n\in \BZ_{\ge 0}$.

{\bf Proof.} One direction is trivial since the natural map $Sym_A^n(X)\to
Sym_A^n(X\oplus Y)$ is injective for any $X$ and $Y$.  

Let us prove another direction.
We are looking for an $A-$algebra $B$ and two maps $\alpha : \be_B
\to M_B$ and $\beta :M_B\to \be_B$ such that $\beta \alpha =\Id$. Now idea is
very simple: let us try to find universal algebra with such properties. For
this we translate our conditions into the language of algebra $B$. First the
map of $B-$modules $\beta :M_B\to \be_B$ is the same as the map of $A-$modules
$M\to \be_B$ or $v:M\to B$ which is equivalent to the map of $A-$algebras
$v_{alg}: Sym^*_A(M)\to B$. Similarly, to give the map $\alpha :\be_B \to M_B$
is the same as to give the map of $A-$modules $u: M^\vee \to B$ or the map of
$A-$algebras $u_{alg}: Sym^*_A(M^\vee )\to
B$. The equation $\beta \alpha =\Id$ is translated to the following condition:
the map
$$\be_A\stackrel{coev}{\longrightarrow} M\otimes M^\vee \stackrel{v\otimes u}
{\longrightarrow} B\otimes B\stackrel{multiplication}{\longrightarrow} B 
\eqno(*)$$
coincides with the map $A\to B$ coming from the fact that $B$ is $A-$algebra.

Summarizing we can say that universal algebra $B$ can be described by
``generators'' $M\oplus M^\vee$ and ``relation'': the map $(*)$ coincides with
the unit morphism $A\to B$ (so algebra $B$ is a
quotient of algebra $Sym_A^*(M\oplus M^\vee)=Sym_A^*(M)\otimes_A
Sym_A^*(M^\vee)$ by the ``ideal'' generated by the morphism $1-\delta :\be_A\to
Sym^*_A(M\oplus M^\vee)$ where $1:\be_A\to Sym^0_A(M)\otimes Sym^0_A(M^\vee)$ 
is the unit morphism and $\delta :\be_A\to Sym_A^1(M)\otimes Sym_A^1(M^\vee)$ 
is the coevaluation morphism). The only problem now is to show that algebra 
described by such generators and relations is {\em nonzero}. 

For this it is enough to show that 1 does not lie in the ideal generated by
$1-\delta$. Assume converse, that is $1=(1-\delta)x$ (here $x$ is a morphism
$\be_A\to Sym_A^*(M\oplus M^\vee)$). Algebra $Sym_A^*(M)
\otimes_ASym_A^*(M^\vee)$ has natural grading by the group $\BZ \oplus \BZ$;
1 lies in $(0,0)-$graded component and $\delta$ lies in $(1,1)-$graded
component. Decompose $x$ in the sum of graded component $x=x_{0,0}+x_{1,0}+
x_{0,1}+x_{1,1}+\ldots$. Clearly, we can assume that $x_{p,q}=0$ for
$p\ne q$ (since if $x$ is solution of $1=(1-\delta)x$ then $x'=x_{0,0}+x_{1,1}
+\ldots$ is a solution too). Now the equation $1=(1-\delta)x$ is equivalent
to the following graded equations:
$$x_{0,0}=1;\; x_{1,1}-\delta x_{0,0}=0;\; x_{2,2}-\delta x_{1,1}=0; \ldots$$
This means that $x_{p,p}=\delta^p$ and $\delta^n=0$ for large enough $n$
since the sum $x=x_{0,0}+x_{1,1}+\ldots$ is finite. Conversely, if
$\delta^n=0$ then $1=(1-\delta)(1+\delta +\ldots +\delta^{n-1})$.

So the universal algebra $B$ is nontrivial if and only if $\delta^n\ne 0$
for all $n$. Now $\delta^n :\be_A\to Sym_A^n(M)\otimes_ASym_A^n(M^\vee)$ equals
to the coevaluation map $\be_A\to Sym_A^n(M)\otimes_ASym_A^n(M)^\vee$ and is
zero if and only if $Sym^n_A(M)=0$. The Lemma is proved. \sq   

\subsection{Local properties} 
One says that a system of objects and
morphisms has some property {\em locally} if this property holds after some
extension of scalars (here the word ``locally'' refers to the topology fppf ---
fidelement plat de presentation finie). For example two objects $X$ and $Y$
of the category $\cA$ are locally isomorphic if there exists a (nonzero) 
algebra $A$ such that $X\otimes A$ is isomorphic to $Y\otimes A$ as 
$A-$modules.

{\bf Exercise.} Let $G$ be an affine algebraic group. Two objects $X, Y\in
\Rep(G)$ are locally isomorphic if and only if $\dim X=\dim Y$. Hint: 
consider algebra of functions on the affine variety $Isom(X,Y)$ --- open
subset of the vector space $X^*\otimes Y$ consisting of isomorphisms; this
variety has natural $G-$action, so algebra of functions has natural structure
of algebra in $\Rep(G)$.  

{\bf Example.} Any short exact sequence in the category $\cA$ locally splits.
Indeed first we reduce ourselves to exact sequence of the form $0\to X\to Y
\stackrel{b}{\to}
\be \to 0$ by the standard argument: splitting of the sequence $0\to M\to N\to
P\to 0$ is equivalent to the splitting of the sequence $0\to M\otimes P^\vee
\to E\to \be \to 0$ where $E$ is the preimage of $\be \in P\otimes P^\vee$ 
under the map $N\otimes P^\vee \to P\otimes P^\vee$. Now one proceeds similarly
(but easier) to the Key Lemma: it is enough to show that algebra 
$Sym^*(Y^\vee)/(b^t-1)$ (where $b^t:\be \to Y^\vee$ is the morphism dual to
$b:Y\to \be$) is nonzero; as before this is equivalent to nonnilpotency of
$b^t$ which is obvious.

Now suppose that category $\cA$ contains an object $\bo$ such that 
$\bo \otimes \bo$ is isomorphic to $\be$ and the commutativity morphism
$\bo \otimes \bo \to \bo \otimes \bo$ is the multiplication by $(-1)$.
Such an object allows to define tensor functor $F:s\Ve \to \cA$ by the
formula $F(V)=V^0\otimes \be \oplus V^1\otimes \bo$ which is an equivalence 
of the category $s\Ve$ and the subcategory $<\be,\bo>$ of $\cA$ consisting
of direct sums of $\be$ and $\bo$.

{\bf Proposition.} For an object $X\in \cA$ the following conditions are
equivalent:

(i) There exist $p$ and $q$ such that $X$ is locally isomorphic to $\be^p\oplus
\bo^q$.

(ii) There exists the Schur functor $S_\lambda$ such that $S_\lambda (X)=0$.

{\bf Proof.} (i)$\Rightarrow$(ii) is trivial. Let us prove that (ii)
$\Rightarrow$ (i). Assume that after some extension of scalars we have
$X_A=\be_A^r\oplus \bo_A^s\oplus S$ for some $A-$module $S$. Consider three
cases:

(a) $Sym_A^n(S)\ne 0$ for all $n$. Then using Key Lemma we can find $A-$algebra
$B$ such that $S_B$ has $\be_B$ as a direct summand and we get decomposition
$X_B=\be_B^{r+1}\oplus \bo_B^s\oplus S'$.

(b) $Sym_A^n(\bo \otimes S)=\bo^{\otimes n}\otimes \wedge_A^n(S)\ne 0$ for all 
$n$. Then again using Key Lemma one finds $A-$algebra $B$ such that
$X_B=\be_B^r\oplus \bo_B^{s+1}\oplus S'$.

(c) Neither (a) nor (b) is true, that is there are $n$ and $m$ such that
$Sym_A^{n+1}(S)=\wedge_A^{m+1}(S)=0$. Let $k$ be any integer greater than
$mn$ and let $\lambda$ be a partition of $k$. Then it is easy to see that
$S_\lambda(S)=0$ (since partition $\lambda$ contains either row of length
greater than $n$ or column of length greater than $m$; in other words any
representation of $S_k$ contains either trivial representation of $S_n$ or
sign representation of $S_m$). Hence $S^{\otimes k}=
\bigoplus_\lambda V_\lambda \otimes S_\lambda(S)=0$. Hence $S=0$ ($S$ is a
direct summand of rigid module, so is rigid; for a rigid module $S$ the
equality $S\otimes_AM=0$ implies $M=0$).

Now we apply iteratively constructions (a) and (b) beginning from the case
$A=\be, r=s=0, S=X$. If this process never stops then $X$ locally has a 
direct summand $\be^p\oplus \bo^q$ with $p+q$ arbitrarily large. But this 
contradicts to the condition $S_\lambda (X)=0$ for some $\lambda$ (see
Exercise 2.2). So at some moment we arrive at (c) and get that $X$ is locally
isomorphic to $\be^p\oplus \bo^q$. \sq

\subsection{Super fiber functor over big ring} Now we can prove that the
super fiber functor exists over sufficiently big (infinitely generated) 
$k-$algebra.

{\bf Proposition.} Assume that any object of the category $\cA$ is annihilated
by some (depending on object) Schur functor. Then there exists a nonzero
supercommutative $k-$algebra $R$ and the $R-$fiber functor of the category 
$\cA$.

{\bf Proof.} We can (and will) assume that the category $\cA$ contains object 
$\bo$ with the properties above (otherwise consider category $\cA_1=\cA
\boxtimes s\Ve$, since $\cA_1$ contains $\cA$ the existence of fiber functor 
for $\cA_1$ implies existence of fiber functor for $\cA$). 

By 2.2 we know that for any object $X$ of $\cA$ there exists an
algebra $B$ such that $X_B=\be_B^r\oplus \bo_B^s$ and for any short exact
sequence in $\cA$ there exists algebra $B$ such that this sequence splits
after extension of scalars to $B$. Let $A$ be the tensor product of all such
algebras (so we need to consider the tensor product of infinitely many 
algebras; this is just inductive limit of tensor products of finitely many
algebras). Then after extension of scalars to $A$ any short exact sequence
in $\cA$ splits and for any $X\in \cA$ one has $X_A=\be_A^r\oplus \bo^s_A$ (in
particular this means that a superdimension of objects of $\cA$ is well
defined).

Consider the functor $\rho :\cA \to s\Ve$ defined by $\rho (X)_0=\Hom(\be,X)$
and $\rho(X)_1=\Hom(\bo,X)$ (this functor is not exact in general). Clearly
$\rho(A)$ is a supercommutative algebra and if $M$ is $A-$module then 
$\rho(M)$ is $\rho(A)-$module (since $<\be,\bo>$ is tensor subcategory of 
$\cA$).
Moreover for two $A-$modules $M, N$ the canonical morphism $can: \rho(M)
\otimes_{\rho(A)}\rho(N)\to \rho(M\otimes_AN)$ is defined. Note that if
$M$ has the form $A\otimes M_0$ with $M_0\in <\be,\bo>$ then $\rho(M)=\rho(A)
\otimes M_0$ (so for any $X\in \cA$ the $\rho(A)-$module $\rho(X_A)$ is free);
if $N$ also has the form $A\otimes N_0$ with $N_0\in <\be,\bo>$ then $M
\otimes_A
N=A\otimes (M_0\otimes N_0)$ and the morphism $can: \rho(M)\otimes_{\rho(A)}
\rho(N)\to \rho(M\otimes_AN)$ is isomorphism. 

Now set $R:=\rho(A)$ and define the functor $\omega: \cA \to R-mod$ by the
formula $\omega (X):=\rho(X_A)=\rho(X\otimes A)$. The remarks above show
that this functor has a natural structure of tensor functor. Moreover since
any short exact sequence in $\cA$ splits after extension of scalars to $A$
this functor is exact. The Proposition is proved. \sq

\subsection{From $R-$fiber functor to a super fiber functor} In this section 
we explain how to get a super fiber functor $F$ from an $R-$fiber functor 
$\omega$. Very roughly the idea is the following: let $X$ be a 
$\otimes-$generator of
$\cA$. Since the superdimension of $X$ is well defined (see previous section)
the superspace $F(X)$ is uniquely determined. Now one needs to define
sufficiently many maps between spaces $F(X)^{\otimes n}$ to ensure that
$F(Y)$ is defined for every $Y\in \cA$ and axioms of super fiber functor are
satisfied (recall that any object of $\cA$ is a subquotient of a direct sum
of objects of the form $X^{\otimes n}$). This is a problem with countably
many variables and countably many equations and it has solution in some algebra
$R$; this implies that this problem has a solution over a field $k$ (see 
\cite{D} for precise statements and proofs). 
 
We restrict ourselves to the case when category $\cA$ is semisimple and has 
only finitely many simple objects $X_1, \ldots ,X_n$. Again the super spaces
$F(X_i)$ are uniquely determined and the only problem is to define tensor
structure on the functor $F$. For this let us introduce variables which
describe all possible isomorphisms $F(X_i)\otimes F(X_j)\to F(X_i\otimes X_j)$
(where in the RHS we use decomposition of $X_i\otimes X_j$ into the sum of
simple objects); these variables (there are only finitely many of them) should
satisfy finitely many equations which express the fact that isomorphisms 
above commute with associativity, commutativity and unit constraints. From
the existence of $R-$fiber functor we know that this problem has a solution
with values in $R$, so it has solution with values in some finitely
generated algebra (since there are only finitely many variables), hence it
has solution with values in $k$ (since by Hilbert Nullstellensatz any
finitely generated algebra over $k$ admits a homomorphism to $k$). Hence
the super fiber functor exists.


\begin{thebibliography}{999}

\bibitem{D} P.~Deligne, {\em Cat\'egories tensorielles}, Moscow Math. Journal
{\bf 2} (2002) no. 2, 227-248.

\bibitem{D1} P.~Deligne, J.~S.~Milne, {\em Tannakian categories}, in ``Hodge
cycles, motives and Shimura varieties, LNM 900, 101-228.

\bibitem{De2} P.~Deligne, {\em La s\'erie exceptionnelle de groupes de Lie},
C.~R.~Acad.~Sci.~Paris S\'er. I Math. {\bf 322} (1996), no. 4, 321-326.

\bibitem{EG1} P.~Etingof, S.~Gelaki, {\em Some properties of finite-dimensional
semi-simple Hopf algebras}, Math. Res. Lett. {\bf 5} (1998), 191-197.

\bibitem{EG2} P.~Etingof, S.~Gelaki, {\em The classification of 
finite-dimensional triangular Hopf algebras over an algebraically closed
field of characteristic 0}, math.QA/0202258.

\bibitem{Mc} S.~Mac Lane, {\em Category theory for the working mathematician},
Springer-Verlag, New York.

\end{thebibliography}
\end{document}